\documentclass[conference]{IEEEtran}
\usepackage{cite}
\usepackage{amsmath,amssymb,amsfonts}
\usepackage{algorithm}
\usepackage[noend]{algpseudocode}
\algrenewcommand\algorithmicrequire{%
    \textbf{\sffamily Input:}}
\algrenewcommand\algorithmicensure{%
    \textbf{\sffamily Output:}}
\algrenewcommand\algorithmiccomment[1]{%
    \hfill\textcolor{gray}{\scriptsize$\triangleright$~#1}}
\newcommand{\CPU}{\textbf{\sffamily CPU}}
\newcommand{\GPU}{\textbf{\sffamily GPU}}
\usepackage{graphicx}
\usepackage{booktabs}
\usepackage{multirow}
\usepackage{xcolor}
\usepackage{tikz}
\usetikzlibrary{arrows.meta,calc,positioning,shapes.geometric}
\usepackage{url}
\usepackage{balance}
\usepackage[caption=false,font=footnotesize]{subfig}
\usepackage{xcolor}
\usepackage[
    colorlinks=true,
    linkcolor=red,     % Algorithms, equations, figures, tables
    citecolor=blue,    % References
    urlcolor=blue      % URLs
]{hyperref}
\usepackage[nameinlink,capitalise,noabbrev]{cleveref}

\IEEEoverridecommandlockouts

\begin{document}
\title{Covering Location Optimization at GPU Speed}

\DeclareRobustCommand*{\IEEEauthorrefmark}[1]{%
    \raisebox{0pt}[0pt][0pt]{\textsuperscript{\footnotesize\ensuremath{#1}}}}
\author{\IEEEauthorblockN{Yun-Tian Zhang and Chen Chen} 
\IEEEauthorblockA{National Key Lab of Autonomous Intelligent Unmanned Systems, Beijing Institute of Technology}
\IEEEauthorblockA{Email: \{zhangyuntian, xiaofan\}@bit.edu.cn }
}

\maketitle

\begin{abstract}
Large-scale covering location optimization often involves NP-hard combinatorial problems that pose substantial computational challenges in urban facility planning. This paper proposes a GPU-accelerated parallel search framework for the maximal covering location problem, combining device-resident multi-trajectory exploration with lightweight CPU-side coordination. Exploiting the structure of one-swap moves, exact incremental gains are evaluated through sparse coverage representations and tiled CUDA kernels. A hierarchical reduction identifies the best exchange without materializing the full swap-gain matrix, while GPU-resident destroy-and-repair and elite-guided restarts balance intensification and diversification. Computational experiments on nine real-world urban configurations demonstrate consistent improvements over GRASP, achieving a mean relative gap of 0.029\% against published best-known solutions (BKS) and reaching BKS in 50 of 90 independent runs. On the largest urban benchmark, Kings, with 7,730 candidate locations and a budget of 50 facilities, the proposed method reaches BKS in all ten runs using a consumer-grade NVIDIA GTX 1660 Ti GPU, requiring only 45.2 seconds on average to obtain its final incumbent, compared with 99.5 seconds reported for AILS-PR in a high-performance computing cluster environment. The source code and benchmark reproduction scripts are publicly available at \url{https://github.com/Yuntian-Zhang/GPU4CO}.
\end{abstract}

\section{Introduction}

Facility location decisions arise in emergency response, public services, and urban infrastructure, where limited resources must serve geographically distributed demand~\cite{murray2016,church2018,saldanha2022}. The maximal covering location problem (MCLP), introduced by Church and ReVelle~\cite{church1974}, selects a prescribed number of facilities to maximize demand-weighted coverage within a service radius. Although its formulation is compact, selecting $p$ out of $n$ candidate sites gives $\binom{n}{p}$ possible configurations. The resulting combinatorial complexity~\cite{bordeaux2014} motivates the search for high-quality solutions under practical computational budgets.

A rich literature has developed exact and heuristic methods for the MCLP. On the exact side, Downs and Camm~\cite{downs1996} combined Lagrangian bounds with branch-and-bound, while Cordeau \emph{et al.}~\cite{cordeau2019} introduced an efficient branch-and-Benders-cut method for large-scale covering location problems. More recently, Chen \emph{et al.}~\cite{chen2023} strengthened mixed-integer optimization through specialized presolving. Zhang \emph{et al.}~\cite{zhang2024integer} studied reformulation-linearization techniques for a probabilistic angular variant of MCLP. These advances substantially extend the capabilities of exact algorithms and general-purpose solvers. Obtaining certificates of optimality can nevertheless be costly when a location decision must be made within a short time budget. Their efficiency also remains sensitive to instance structure and size: in particular, the Benders formulation of~\cite{cordeau2019} is especially suited to settings where demand points far outnumber candidate sites, whereas our settings involve large numbers of both demand points and candidate facilities.

Heuristic research provides a complementary line of development. Early greedy and interchange procedures~\cite{church1974} were followed by Lagrangian-based search~\cite{galvao1996}, GRASP~\cite{resende1998,resende2016,laguna2025}, heuristic concentration~\cite{revelle2008}, and decomposition heuristics~\cite{senne2010}. More recent approaches include intelligent-guided adaptive search~\cite{maximo2017} and hybrid adaptive iterated local search with path relinking (AILS-PR)~\cite{maximo2025}. While these methods offer increasingly effective search strategies, traditional implementations largely depend on CPU-oriented computation, as do conventional mixed-integer programming solvers. In particular, evaluating large neighborhoods requires extensive computation over overlapping coverage relationships, even when the underlying search strategy is sophisticated. \emph{Scaling} to large neighborhoods therefore raises a computational question beyond the choice of search strategy itself: how can the many candidate moves around an incumbent be evaluated efficiently?

This question naturally motivates GPU-oriented search. Recent frameworks, including EvoX~\cite{huang2025evox} and cuGenOpt~\cite{liu2026cugenopt}, demonstrate the potential of GPU acceleration for evolutionary computation and general-purpose metaheuristics. Rather than pursuing a problem-agnostic framework, we exploit the specific structure of weighted coverage. A one-swap neighborhood contains $p(n-p)$ exchanges; for Kings with $n=7{,}730$ and $p=50$, there are $384{,}000$ moves per incumbent. Sparse coverage representations and incremental gains make these exchanges amenable to parallel evaluation, with many candidates processed by CUDA thread blocks and their best moves selected through device-side reductions. However, speedup requires that trajectory states remain on the device: repeated host--device transfers and synchronization can otherwise erode the benefits of massively parallel neighborhood scans.

We propose a GPU-accelerated parallel search framework combining device-resident neighborhood exploration with lightweight CPU-side coordination. Sparse incremental evaluation and tiled CUDA kernels exploit weighted coverage structure, while destroy-and-repair and elite-guided restarts balance intensification and diversification with limited host--device communication. Experiments on nine urban MCLP configurations demonstrate competitive performance against GRASP and COPT~\cite{copt}, with comparisons to published AILS-PR results~\cite{maximo2025}. Implemented using Numba-CUDA~\cite{cuda,numba}, the method reaches best-known solutions in 50 of 90 runs with a mean gap of $0.029\%$. Notably, a consumer-grade NVIDIA GTX 1660 Ti reaches the best-known solution (BKS) in all ten runs on the largest benchmark, demonstrating the potential of accessible GPU hardware for large-scale combinatorial optimization in practical applications.

\section{Problem Formulation and Search Representation}

\subsection{Problem Formulation}

Let $I$ be the demand point set and $J$ the candidate facility
set, with $|I|=m$ and $|J|=n$. Each demand point $i\in I$
has nonnegative weight $w_i$. Define the binary incidence
$a_{ij}=1$ if placing a facility at $j$ covers demand $i$,
and $a_{ij}=0$ otherwise. The standard binary formulation~\cite{church1974} is
\begin{align}
 \max_{\boldsymbol{x,y}}\quad
 & \sum_{i\in I}w_i y_i, \label{eq:obj}\\
 \mathrm{s.t.}\quad
 & y_i\le \sum_{j\in J}a_{ij}x_j,
 \quad i\in I, \label{eq:cover}\\
 & \sum_{j\in J}x_j=p, \label{eq:card}\\
 & x_j\in\{0,1\},\quad y_i\in\{0,1\}.
 \label{eq:binary}
\end{align}

Because all weights are nonnegative, a selected facility set
$S=\{j:x_j=1\}$ has objective
\begin{equation}
 F(S)=\sum_{i\in I}w_i\,
 \mathbf{1}\!\left\{
 \sum_{j\in S}a_{ij}\ge1
 \right\},\qquad |S|=p.
 \label{eq:setobj}
\end{equation}

The search operates directly on $S$ and its coverage state.
A facility-oriented compressed sparse row (CSR) representation
stores the covered-demand indices of each candidate site,
rather than the dense incidence matrix. For the current
selection, the coverage count
$c_i=\sum_{j\in S}a_{ij}$ records how many selected facilities
cover demand point $i$. The selected sites, membership
indicators, and coverage counts form the state maintained
by each search trajectory.

\subsection{Exact Gain of a One-Swap Move}

Consider removing $r\in S$ and adding $s\notin S$, giving
$S'=(S\setminus\{r\})\cup\{s\}$.
The objective change can be decomposed into three quantities.

The \emph{removal loss} is
$\mathcal{L}_r=
\sum_{i:\,a_{ir}=1,\,c_i=1}w_i$,
representing demand covered exclusively by $r$.
The \emph{addition gain} is
$\mathcal{G}_s=
\sum_{i:\,a_{is}=1,\,c_i=0}w_i$,
representing previously uncovered demand reached by $s$.
Finally, the \emph{overlap correction}
$\mathcal{H}_{rs}=
\sum_{i:\,a_{ir}=a_{is}=1,\,c_i=1}w_i$
restores the removal loss for demand that $s$ continues
to cover.

The exact one-swap gain is therefore
\begin{equation}
 \Delta(r,s)=F(S')-F(S)
 =\mathcal{G}_s+\mathcal{H}_{rs}-\mathcal{L}_r.
 \label{eq:swap}
\end{equation}

The removal losses can be precomputed and reused across
candidate additions. Moreover, evaluating $s$ requires
visiting only its covered-demand indices, making the gain
calculation well suited to sparse incremental evaluation.

\noindent\textbf{Illustrative example.}
Consider four demand points with weights $(5,3,4,2)$ and
three candidate facilities:
$N(A)=\{1,2\}$,
$N(B)=\{2,3\}$, and
$N(C)=\{1,4\}$,
where $N(j)$ denotes the demand covered by facility $j$.
For $p=2$ and $S=\{A,B\}$, the coverage counts are
$c=(1,2,1,0)$ and $F(S)=12$.
When adding $C$, we obtain
$\mathcal{G}_C=2$,
$\mathcal{L}_A=5$,
$\mathcal{L}_B=4$,
$\mathcal{H}_{AC}=5$, and
$\mathcal{H}_{BC}=0$.
Consequently,
$\Delta(A,C)=2+5-5=2$, whereas
$\Delta(B,C)=2+0-4=-2$.
The improving swap $A\rightarrow C$ gives
$S'=\{B,C\}$ and $F(S')=14$, without recomputing
the full coverage objective.

A best-improvement step evaluates all feasible one-swaps
and applies the largest positive gain. Although this
neighborhood admits exact incremental evaluation, a
one-swap local optimum need not be globally optimal.
Rather than exhaustively enumerating higher-order swaps,
the proposed search employs destroy-and-repair to explore
beyond such local optima.

\section{GPU-Accelerated Parallel Search}
\label{sec:gpu-search}

\subsection{Device-Resident Multi-Trajectory Search}

The proposed method combines massively parallel neighborhood
evaluation on the GPU with lightweight search coordination
on the CPU. It maintains $K=32$ concurrent trajectories,
each with a selected set $S_b$, coverage counts $c_{bi}$,
objective $F(S_b)$, and an independent best-so-far solution.
The sparse incidence data and trajectory states remain
resident on the GPU throughout each search pool.

\begin{figure*}[htbp]
    \centering
    \includegraphics[width=0.95\textwidth]{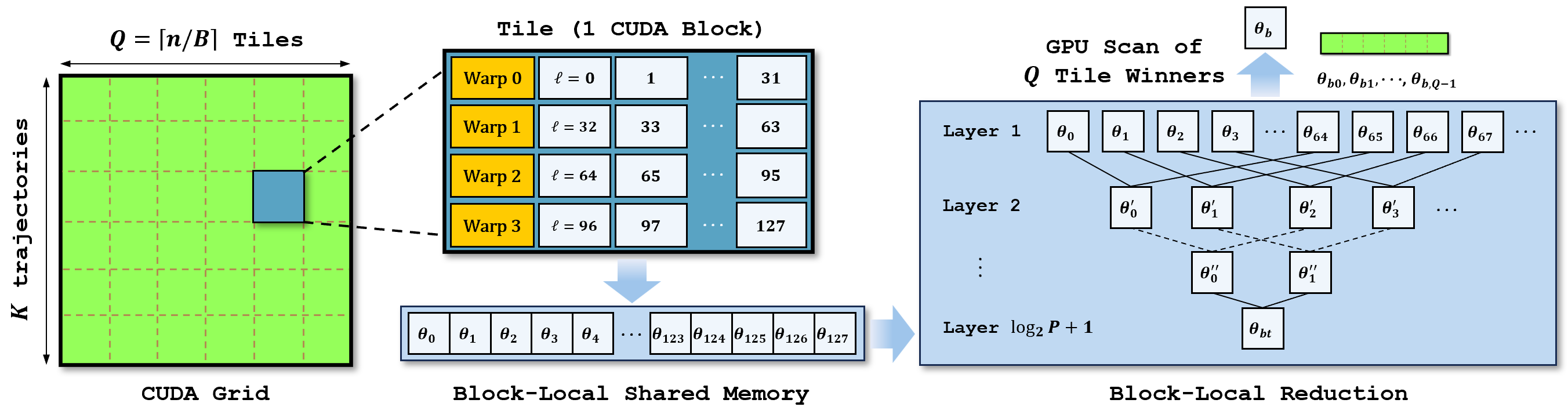}
    \caption{CUDA mapping and hierarchical best-swap selection.
    Each trajectory--tile pair is assigned to a 128-thread
    block. Thread-local best swaps are written to block-local
    shared memory and reduced to one tile winner. A subsequent
    GPU scan selects the best swap across tiles for each
    trajectory.}
    \label{fig:gpu}
\end{figure*}

Figure~\ref{fig:gpu} illustrates the CUDA execution
structure. Parallelism is exposed at two levels:
independent search trajectories and candidate-site tiles
within each trajectory. The GPU performs exact one-swap
evaluation, best-move selection, and incremental state
updates; the CPU manages search pools, elite solutions,
and restart initialization. This separation reduces
host--device transfers during the computationally intensive
local-search iterations.

\subsection{GPU-Tiled Best-Improvement Evaluation}

For each trajectory $b$, the GPU first identifies the unique
covering facility of demand points with $c_{bi}=1$ and
precomputes removal losses $\mathcal{L}_{br}$.
Each candidate addition $s$ then requires only a CSR traversal
of $N(s)=\{i\in I:a_{is}=1\}$ to accumulate its addition gain
$\mathcal{G}_{bs}$ and overlap corrections $\mathcal{H}_{brs}$.
The best removal position is identified from the $p$
possibilities using~\eqref{eq:swap}.

To parallelize this computation, candidate-site indices
are partitioned into tiles of $B=256$ sites. The two-dimensional
CUDA grid is indexed by trajectory $b$ and tile $t$,
with one CUDA block assigned to each pair $(b,t)$,
as shown in Fig.~\ref{fig:gpu}.
Each block contains $P=128$ threads, organized into
four warps of 32 threads. A thread evaluates up to two
candidate additions and retains its best swap triple
$\theta_\ell=(\Delta_\ell,r_\ell,s_\ell)$.

The thread-local triples are written to block-local
shared memory. A tree reduction compares them using
successively halved strides,
$64,32,\ldots,1$, requiring seven comparison stages
to reduce 128 candidates to one tile winner
$\theta_{bt}$. This block-level reduction has linear
total work and logarithmic parallel depth in shared memory.

The tile winners are stored in GPU global memory.
A subsequent GPU kernel assigns one thread to each
trajectory to scan its $Q=\lceil n/B\rceil$ tile winners,
yielding the trajectory-best swap $\theta_b$.
Only a positive-gain swap is applied, followed by
incremental updates of the selected facilities and
coverage counts. Thus, the complete one-swap neighborhood
is evaluated without materializing a dense gain matrix.

\begin{algorithm}[htbp]
\caption{\textsf{GPU Tiled Best-Swap Evaluation}}
\label{alg:tile}
\footnotesize
\begin{algorithmic}[1]
\Require Resident states $S_b,c_b$; CSR lists $N(s)$; losses $\mathcal{L}_{br}$; $B=256$, $P=128$
\Ensure Best swap triple $\theta_b=(\Delta_b,r_b,s_b)$ for each trajectory
\ForAll{$(b,t)$ \textbf{in parallel}: trajectory $b$, tile $t$}
  \ForAll{thread $\ell\in\{0,\ldots,P-1\}$ \textbf{in parallel}}
    \State $\theta_\ell\gets(-\infty,-1,-1)$
    \For{$s\in\{tB+\ell,\ tB+\ell+P\}\cap J$}
      \If{$s\notin S_b$}
        \State $\mathcal{G}_{bs}\gets0$; $\mathcal{H}_{brs}\gets0$ for each $r\in S_b$
        \ForAll{$i\in N(s)$ \textbf{via CSR}}
          \State $\mathcal{G}_{bs}\gets\mathcal{G}_{bs}+w_i$ if $c_{bi}=0$
          \State $\mathcal{H}_{brs}\gets\mathcal{H}_{brs}+w_i$ if $c_{bi}=1$ and $r$ is its unique owner
        \EndFor
        \State Update $\theta_\ell$ with $\max_{r\in S_b}(\mathcal{G}_{bs}+\mathcal{H}_{brs}-\mathcal{L}_{br},r,s)$
      \EndIf
    \EndFor
  \EndFor
  \State $\theta_{bt}\gets\Call{BlockArgmax}{\theta_0,\ldots,\theta_{P-1}}$ \Comment{shared memory}
\EndFor
\ForAll{trajectory $b$ \textbf{in parallel}}
  \State $\theta_b\gets\Call{ScanArgmax}{\theta_{b0},\ldots,\theta_{b,Q-1}}$ \Comment{on GPU}
\EndFor
\end{algorithmic}
\end{algorithm}

For Kings ($n=7{,}730$), each trajectory has $Q=31$ tiles.
The resulting launch grid contains
$32\times31=992$ CUDA blocks, each with 128 threads.
Algorithm~\ref{alg:tile} summarizes the device-side
evaluation and hierarchical best-move selection.

\subsection{Escaping Local Optima and Elite-Guided Restarts}

When a trajectory reaches a one-swap local optimum,
a GPU-resident destroy-and-repair operator perturbs its
current selection. It removes $k=4$ facilities, temporarily
excludes the removed sites during repair, and restores
cardinality by inserting four facilities sequentially.
Each insertion maximizes the newly covered demand weight,
using the same tiled parallel evaluation pattern.
Coverage counts are updated after each insertion.

The perturbed solution may have a lower objective,
allowing the trajectory to explore a different search
basin. Its best-so-far solution is retained separately.
A trajectory becomes inactive after exhausting at most
20 perturbations.

At pool boundaries, the CPU collects and validates
the trajectory-best solutions and maintains an elite
archive $\mathcal E$ of at most 12 high-quality,
diverse selections. Diversity is measured by the
exchange distance
\begin{equation}
 d(S_1,S_2)=p-|S_1\cap S_2|,
 \label{eq:dist}
\end{equation}
with a minimum separation target of 10 swaps.
Near-duplicate solutions are replaced only when improved,
while sufficiently inferior solutions are discarded.

The archive guides the next pool of 32 trajectories:
16 starts are obtained by elite mutation, eight by
elite recombination, and eight are generated afresh.
Mutation replaces 6, 10, 14, or 20 facilities.
Recombination retains facilities common to two elite
parents and completes the selection through randomized
greedy construction. When too few elites are available,
recombination falls back to mutation.

\begin{algorithm}[htbp]
\caption{\textsf{Batch-Synchronized CPU--GPU Search}}
\label{alg:search}
\footnotesize
\begin{algorithmic}[1]
\Require Sparse incidence $A$, demand weights $w$, facilities $p$, time limit $T$
\Ensure Best observed facility set $S^*$ and $F(S^*)$
\State \CPU: Build weighted greedy anchor; set $\mathcal E\gets\varnothing$
\State \CPU: Allocate shared instance data and $K=32$ device states
\While{time budget remains}
  \State \CPU: Generate $K$ starts (fresh / elite mutation / recombination)
  \State \GPU: Initialize $S_b,c_b,F(S_b)$ and trajectory-best records
  \While{some trajectory is active and time remains}
    \State \GPU: Compute unique owners and removal losses $\mathcal{L}_{br}$
    \State \GPU: $\theta_b\gets\Call{TiledBestSwap}{S_b,c_b,\mathcal{L}_b}$ \Comment{Alg.~\ref{alg:tile}}
    \ForAll{active trajectories $b$ \textbf{in parallel}}
      \If{best gain $\Delta_b>0$}
        \State Apply swap; update $c_b$, $F(S_b)$, and best record
      \ElsIf{fewer than 20 perturbations have been used}
        \State Destroy four sites; greedily repair via GPU tiles
      \Else
        \State Mark trajectory $b$ inactive
      \EndIf
    \EndFor
    \State \CPU: Synchronize status every five rounds
  \EndWhile
  \State \CPU: Collect and validate trajectory-best solutions
  \State \CPU: Update $\mathcal E$; retain the best overall $S^*$
\EndWhile
\State \Return $S^*, F(S^*)$
\end{algorithmic}
\end{algorithm}

% Requires: booktabs, multirow, graphicx
\begin{table*}[htbp]
\centering
\caption{Solution quality and computational time on nine real-world MCLP configurations.}
\label{tab:results}
\renewcommand{\arraystretch}{1.35}
\setlength{\tabcolsep}{1.5pt}
\scriptsize
\resizebox{\textwidth}{!}{%
\begin{tabular}{llc ccc ccc cccc}
\toprule
\multirow{2}{*}{Instance}
& \multirow{2}{*}{$p$}
& \multirow{2}{*}{BKS}
& \multicolumn{3}{c}{Proposed GPU Search}
& \multicolumn{3}{c}{GRASP}
& \multicolumn{4}{c}{COPT} \\
\cmidrule(lr){4-6}
\cmidrule(lr){7-9}
\cmidrule(lr){10-13}
& & &
Gap (\%) $\mu\pm\sigma$ &
BKS hits &
$t_{\mathrm{best}}$ (s) &
Gap (\%) $\mu\pm\sigma$ &
BKS hits &
$t_{\mathrm{best}}$ (s) &
BKS gap (\%) &
Solve time (s) &
Upper bound &
MIP gap (\%) \\
\midrule
\multirow{5}{*}{\shortstack[l]{Kings\\(7,730 nodes)}}
& 50 & 2,022,077
& $0.000\pm0.000$ & 10/10 & 45.212
& $0.530\pm0.083$ & 0/10 & 202.000
& 0.003 & 3600.000 & 2,027,241.645 & 0.258 \\
& 60 & 2,227,025
& $0.034\pm0.026$ & 1/10 & 116.415
& $0.770\pm0.080$ & 0/10 & 172.069
& 1.318 & 3600.000 & 2,254,716.720 & 2.596 \\
& 70 & 2,380,660
& $0.023\pm0.026$ & 4/10 & 90.927
& $1.019\pm0.119$ & 0/10 & 115.481
& 10.850 & 3600.000 & 2,410,231.506 & 13.564 \\
& 80 & 2,467,701
& $0.095\pm0.056$ & 0/10 & 101.409
& $1.066\pm0.051$ & 0/10 & 178.738
& 31.340 & 3600.000 & 2,490,530.739 & 46.993 \\
& 90 & 2,502,669
& $0.077\pm0.024$ & 0/10 & 171.458
& $0.686\pm0.082$ & 0/10 & 166.271
& 22.774 & 3600.000 & 2,504,700.000 & 29.596 \\
\midrule
\multirow{2}{*}{\shortstack[l]{San Francisco\\(5,137 nodes)}}
& 50 & 617,592
& $0.000\pm0.000$ & 10/10 & 10.347
& $0.485\pm0.116$ & 0/10 & 160.003
& 0.000 & 3600.000 & 618,778.345 & 0.192 \\
& 70 & 733,772
& $0.003\pm0.009$ & 9/10 & 121.468
& $0.858\pm0.087$ & 0/10 & 170.339
& 0.215 & 3600.000 & 739,555.288 & 1.005 \\
\midrule
\multirow{2}{*}{\shortstack[l]{Bronx\\(3,839 nodes)}}
& 50 & 1,205,051
& $0.000\pm0.000$ & 10/10 & 67.787
& $0.525\pm0.086$ & 0/10 & 114.852
& 0.000 & 700.304 & 1,205,051.000 & 0.000 \\
& 70 & 1,348,232
& $0.028\pm0.041$ & 6/10 & 33.902
& $0.698\pm0.069$ & 0/10 & 132.023
& 0.000 & 3600.000 & 1,351,483.097 & 0.241 \\
\bottomrule
\end{tabular}}
\begin{minipage}{\textwidth}
\footnotesize
\textit{Remark 1:}
BKS denotes the published best-known feasible objective value.
The reference study~\cite{maximo2025} reports experiments conducted on a
104-node computing cluster, with each node equipped with
two Intel Xeon E5-2680 v2 processors (2.8 GHz, 10 cores each)
and 128 GB of DDR3-1866 RAM.
\par\smallskip
\textit{Remark 2:}
GPU and GRASP results are based on 10 independent runs per
configuration, each subject to a \textbf{300-s time limit}.
COPT was executed once per configuration with a
\textbf{3,600-s time limit}. Reported COPT times are actual
solver times; the Bronx instance with $p=50$ was proven
optimal in 700.304 s. The metric $t_{\mathrm{best}}$
denotes the mean time at which each heuristic first
attains its own final incumbent, excluding setup and
JIT compilation.
\end{minipage}
\end{table*}

This mechanism combines the exploitation of high-quality
search regions with continued exploration through
perturbation, recombination, and fresh starts. By preserving
both solution quality and structural diversity, the elite
archive provides informative starting points without
forcing subsequent trajectories into the same search basin.
The initial pool is constructed from a weighted greedy
anchor and randomized variants, ensuring that the search
remains independent of published best-known facility sets.

\subsection{Overall Coordinated Procedure}

Algorithm~\ref{alg:search} summarizes the complete
CPU--GPU procedure and the interaction between its
local exploration and global restart mechanisms.
Within each search pool, the GPU advances multiple
trajectories through exact best-improvement swaps and
destroy-and-repair operations. The CPU synchronizes
compact incumbent and activity information every five
rounds to monitor progress and the remaining time budget.
At pool boundaries, it collects the trajectory-best
solutions, updates the elite archive, and generates
diversified initializations for the next pool.
This batch-synchronized organization keeps the
computationally intensive search on the GPU while
restricting host--device communication to periodic
status checks and pool transitions. The search terminates
upon reaching the prescribed time budget, excluding
instance preparation and CUDA compilation.

\section{Computational Experiments}

\subsection{Instances and Experimental Protocol}

We evaluate the proposed method on the publicly available MCLP benchmark instances~\cite{maximo2025}, covering three urban regions: Kings, San Francisco, and Bronx. The instances contain 7,730, 5,137, and 3,839 candidate/demand locations, respectively. We use the published service radii of 800 m for Kings and 600 m for the other two regions. Demand weights represent population, and coverage matrices follow the benchmark's spherical-distance convention and radius-boundary treatment. We independently recomputed the objectives of the published reference placements and recovered the reported BKS values.

The nine tested configurations use facility budgets $p\in\{50,60,70,80,90\}$ for Kings, and $p\in\{50,70\}$ for San Francisco and Bronx. Experiments were conducted on a computer equipped with an Intel Core i9-9900 CPU (3.10 GHz), 64 GB RAM, and an NVIDIA GeForce GTX 1660 Ti GPU. The proposed GPU search is implemented in Python using NumPy/SciPy sparse structures and Numba-CUDA, with source code and reproduction scripts publicly available\footnote{\url{https://github.com/Yuntian-Zhang/GPU4CO}}. The GRASP comparator~\cite{resende1998,resende2016,laguna2025} uses randomized greedy construction with restricted candidate-list parameter $\alpha=0.15$, followed by one-swap local improvement.

For each configuration, the proposed GPU search and GRASP were each run ten times with seeds 2027--2036 and a 300-s search limit. Setup and JIT compilation time are excluded from this limit. COPT 8.0~\cite{copt} was run once per configuration with a 3,600-s time limit; the table reports actual elapsed solve time. For Bronx with $p=50$, COPT terminated after 700.304 s upon proving optimality. The published BKS~\cite{maximo2025} values are best-known feasible objectives and are used as reference values throughout the evaluation.

\subsection{Solution Quality on Urban Benchmarks}

Let $z_{\mathrm{BKS}}$ denote the published best-known objective and $\bar z$ the mean final objective over the ten independent runs of a heuristic. We calculate the relative gap as
\begin{equation}
\mathrm{Gap}(\%) =
100\,\frac{z_{\mathrm{BKS}}-\bar z}{z_{\mathrm{BKS}}}.
\label{eq:gap}
\end{equation}
Table~\ref{tab:results} provides a comprehensive
comparison of solution quality and computational
performance across all nine configurations and three
methods. For GPU and GRASP, it reports the mean and
standard deviation of run-level gaps to the BKS,
reflecting both average solution quality and variability
across independent runs. The number of BKS hits
indicates solution consistency, while the mean time
to the final incumbent characterizes the progression
of each heuristic within its search budget.
For COPT, the table presents the incumbent gap to BKS,
elapsed solve time, best upper bound, and solver-reported
MIP gap, allowing feasible solution quality and
optimality certification progress to be assessed
separately.

The proposed GPU search reaches BKS in 50 of 90 runs, whereas GRASP does not reach BKS in any run. The GPU mean gap ranges from 0.000\% to 0.095\%, while GRASP's mean gaps range from 0.485\% to 1.066\%. Averaged over all 90 runs, the final relative gap is 0.029\% for GPU search and 0.737\% for GRASP. The GPU hit counts are 15/50 for Kings, 19/20 for San Francisco, and 16/20 for Bronx. Kings at the larger facility budgets remains the most demanding setting: no run reaches BKS for $p=80$ or $p=90$, with mean gaps of 0.095\% and 0.077\%, respectively. The reported $t_{\mathrm{best}}$ is the mean time to each method's own final incumbent, rather than the time to BKS.

\begin{figure}[htbp]
    \centering
    \includegraphics[width=0.715\columnwidth]{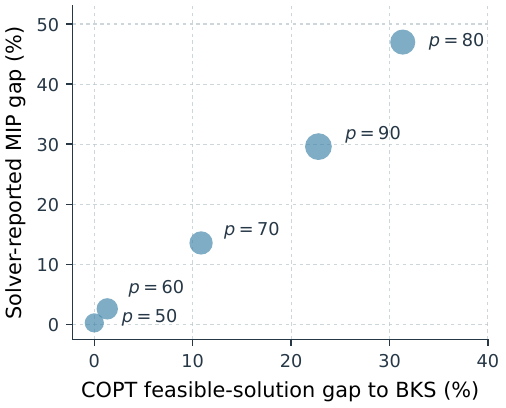}
    \caption{COPT performance on the largest Kings instance ($n=7,730$) under different facility budgets. The horizontal and vertical axes represent the relative gap of the best feasible solution to the BKS and the solver-reported MIP gap, respectively. Bubble sizes increase with the facility budget $p$, as indicated by the labels. Points closer to the origin indicate stronger feasible solutions and tighter optimality bounds.}
    \label{fig:copt}
\end{figure}

COPT reaches a BKS-matching incumbent in three configurations: San Francisco with $p=50$ and Bronx with $p=50$ and $p=70$. It proves optimality only for Bronx with $p=50$, where the incumbent equals the upper bound. The Kings results provide a more demanding comparison across facility budgets. At $p=50$, COPT obtains a solution within 0.003\% of BKS with a 0.258\% MIP gap. The incumbent and MIP gaps increase substantially at $p=70$ and $p=80$, reaching 10.850\% and 13.564\% at $p=70$, and 31.340\% and 46.993\% at $p=80$. Both gaps decrease at $p=90$ but remain substantial. These results distinguish the quality of the incumbent from the difficulty of closing the solver's optimality gap, particularly for larger facility budgets.

\subsection{Exact-Solver Results}

\begin{figure*}[htbp]
    \centering
    \subfloat[Locations with $p=50$ facilities.]{
        \includegraphics[width=0.31\textwidth]{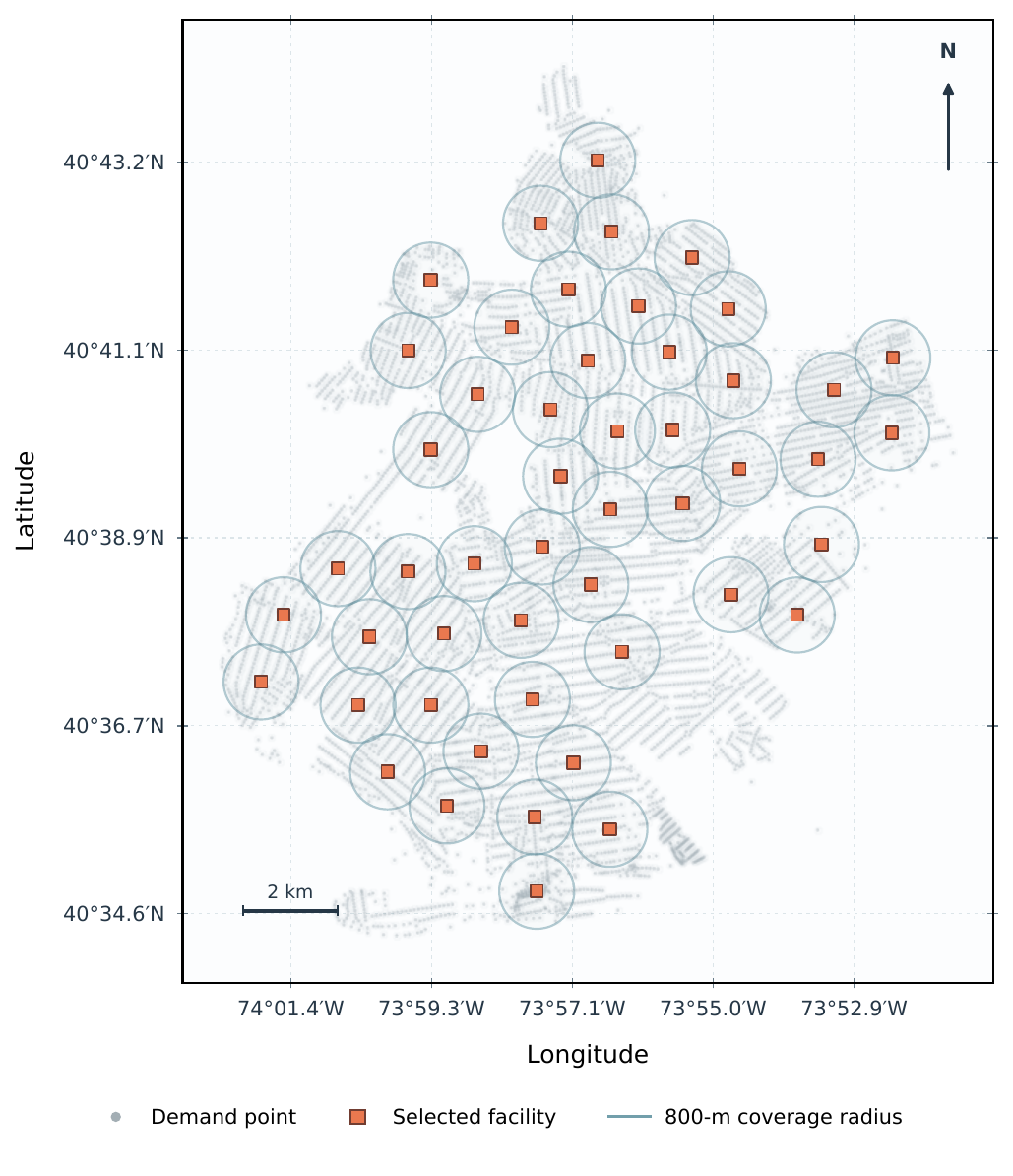}
        \label{fig:kings50}
    }
    \hfill
    \subfloat[Locations with $p=70$ facilities.]{
        \includegraphics[width=0.31\textwidth]{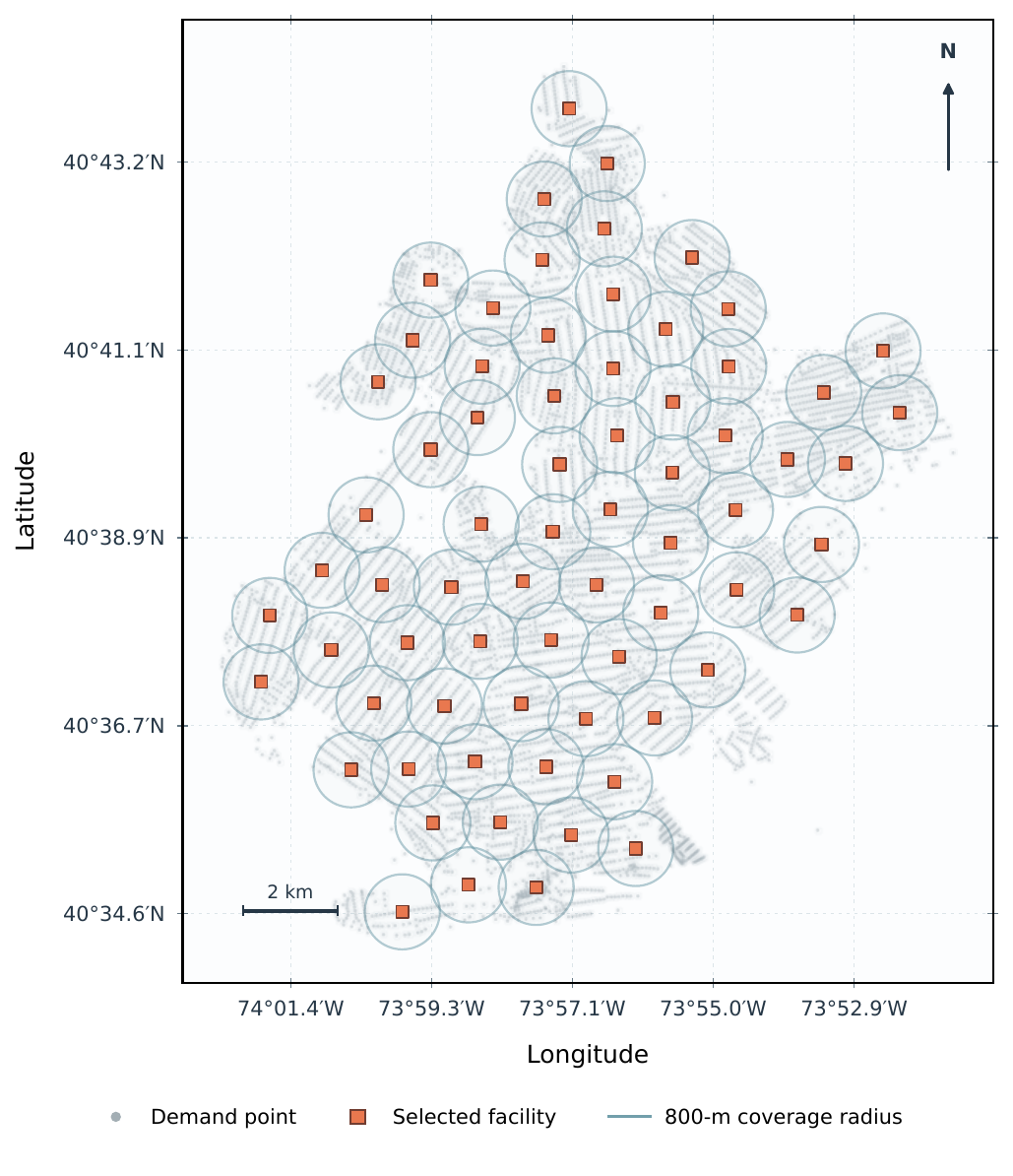}
        \label{fig:kings70}
    }
    \hfill
    \subfloat[Locations with $p=90$ facilities.]{
        \includegraphics[width=0.31\textwidth]{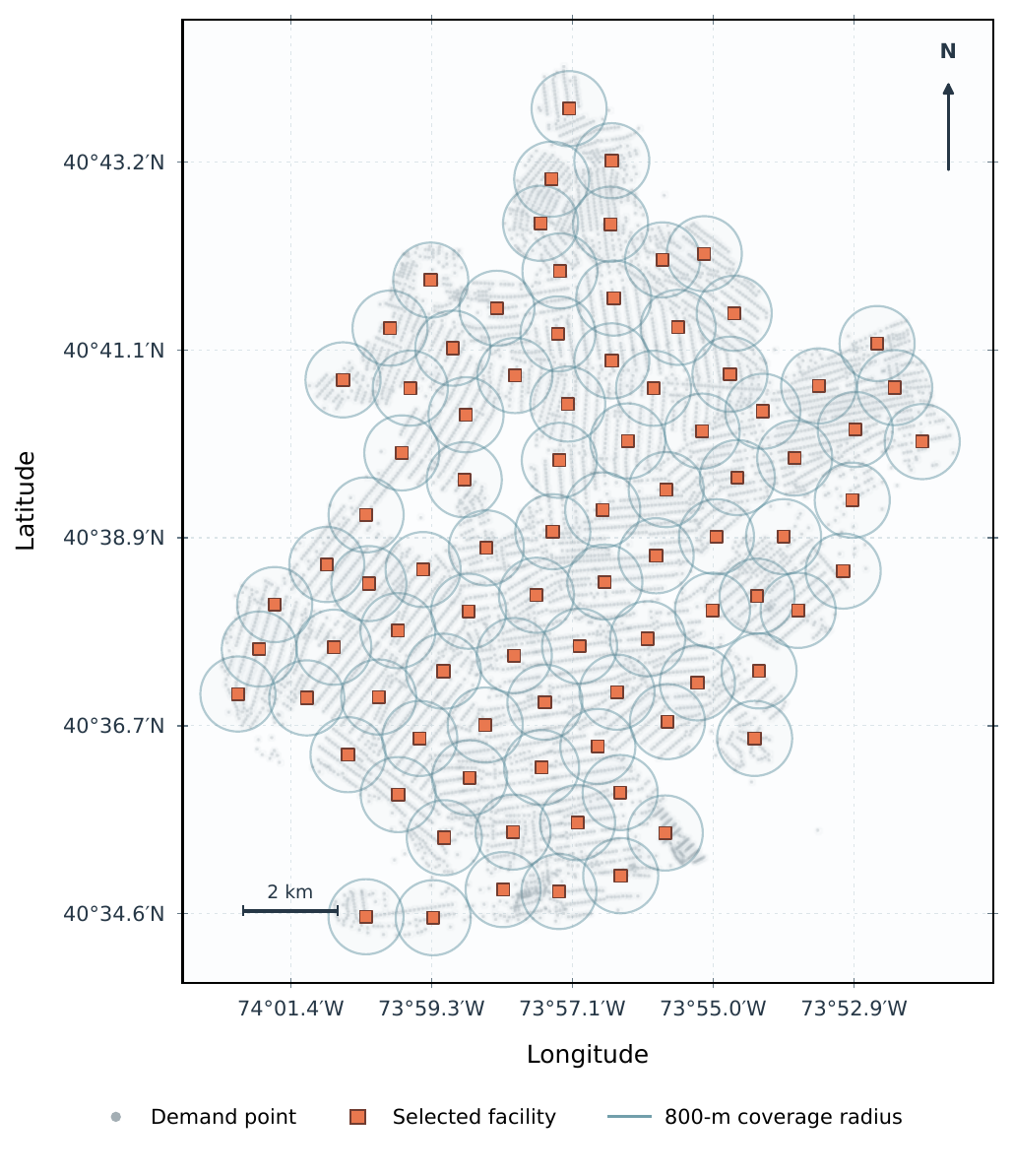}
        \label{fig:kings90}
    }
    \caption{Spatial location patterns of the proposed GPU search on the Kings instance under different facility budgets. Gray circles represent demand points and orange squares denote selected facilities. Circular boundaries illustrate the 800-m coverage radius of each selected facility. All three panels use the same geographic extent, projection, and scale to facilitate spatial comparison.}
    \label{fig:kings-deployments}
\end{figure*}

Beyond the aggregate solution-quality comparisons in
Table~\ref{tab:results}, the exact-solver results
provide further insight into the computational difficulty
of different facility budgets. Fig.~\ref{fig:copt}
focuses on the five Kings configurations, which form
the largest instance family. The two axes distinguish
the quality of the best feasible solution from the
remaining optimality gap.

COPT produces a near-BKS incumbent at $p=50$, with
relatively small gaps also observed at $p=60$.
Both the incumbent and MIP gaps increase substantially
at $p=70$ and reach their largest values at $p=80$.
Although both decrease at $p=90$, they remain considerably
larger than those for smaller budgets. This non-monotonic
pattern suggests that computational difficulty depends
not only on instance size but also on the facility budget
and the resulting coverage interactions. The comparison
further illustrates the value of rapidly obtaining
high-quality feasible solutions when exact optimization
struggles to close the optimality gap.

\subsection{Spatial Location Patterns}

Fig.~\ref{fig:kings-deployments} visualizes the facility locations produced by the proposed GPU search on Kings
for $p=50$, $70$, and $90$. As the facility budget
increases, the selected locations become more densely
distributed across the urban region, with greater spatial
overlap among their coverage areas. The additional facilities provide opportunities to reach previously uncovered demand, while overlapping service regions reflect the spatial concentration of candidate coverage.

These patterns reflect the spatial trade-off inherent
in maximal covering location decisions: additional
facilities can improve the weighted coverage objective,
but their marginal contributions depend on existing
coverage and the geographic distribution of demand.
The maps therefore complement the numerical results
in Table~\ref{tab:results} by illustrating how different
facility budgets translate into actual location
configurations.

\section{Conclusion}

This paper presented a GPU-oriented parallel local search
framework for large-scale covering location optimization.
By combining exact incremental one-swap evaluation,
tiled CUDA neighborhood exploration, and GPU-resident
multi-trajectory search with CPU-side elite-guided
coordination, the proposed method exploits massive
parallelism while limiting host--device communication.
Experiments on nine real-world urban configurations
demonstrated consistently high solution quality,
achieving an average gap of only 0.029\% against published BKS and matching them in 50 of 90 runs.
Notably, on the largest Kings benchmark with $p=50$,
a consumer-grade GTX 1660 Ti reached the best-known
solution in all ten runs, with an average time of
45.2 seconds. These results demonstrate the potential
of accessible GPU hardware for computationally demanding
combinatorial optimization.

This study represents an initial exploration of
GPU-oriented heuristic optimization, leaving substantial
opportunities for further development. One direction is
to investigate richer neighborhoods and more advanced
search mechanisms, while balancing their potential
solution-quality gains against the parallelism, memory
access patterns, and synchronization costs of heterogeneous
CPU--GPU architectures. Beyond algorithmic improvements,
the multi-trajectory design also motivates investigating
\emph{scaling laws} in combinatorial optimization, particularly
how solution quality evolves with GPU computing resources,
parallel search breadth, and computational budgets.
Another promising direction is to extend GPU acceleration
to exact algorithms through batched relaxation solving,
parallel bound computation, and GPU-enhanced
branch-and-bound, potentially combining high-quality
feasible solutions with provable optimality guarantees.

\section*{Acknowledgment}
This research was supported by the National Natural Science Foundation of China under Grant 62673059.

\end{document}